# A New One Parameter Bimodal Skew Logistic Distribution and Its Applications

S. Shah, S. Chakraborty and P. J. Hazarika

Department of Statistics, Dibrugarh University, Dibrugarh, Assam, 786004

**Abstract**

In this paper an attempt is made to develop a new bimodal alpha skew logistic distribution under Balakrishnan (2002) mechanism. Some of its distributional as well as moments properties are studied. Some extensions of this new distribution are also briefly introduced. Finally, suitability of the proposed distribution is checked with the help of data fitting examples of real life.

**Maths Classification:** Primary: 62E15, 60E05; Secondary: 62P05

## 1. Introduction

Although we cannot undermine the applications and value of normal distribution in statistics but there are several random phenomena in nature that we sometimes encounter in datasets having high values of skewness and kurtosis, bimodal and multimodal behaviours in their frequency curve which can't be described by unimodal normal distribution. For example, the asymmetry and bimodality behaviours are seen in many areas of application of data analysis namely, the time between eruptions of certain geysers, the size of worker weaver ants, age of incidence of certain diseases, in age distribution and growth estimates of fish population, in geological studies distribution of the sediments etc. Also, it is observed that the data exhibit many modes as well as asymmetry in the field of demography, insurance, medical sciences, physics, etc. ( for details see Hammel et al. (1983), Dimitrov et al. (1997), Sinha (2012) among others).

To deal with such type of random phenomena, Azzalini (1985) first introduced the path breaking skew normal distribution in his seminal paper entitled "A class of distributions which includes the normal ones. *Scandinavian Journal of Statistics*, 171-178" and its probability density function (pdf) is given by

$$f_Z(z;\lambda) = 2\varphi(z)\Phi(\lambda z); \quad -\infty < z < \infty \qquad (1)$$

where, $\varphi$ and $\Phi$ are, respectively, the pdf and cumulative distribution function (cdf) of standard normal variable $Z$ and $\lambda \in R$, the skewness parameter. Symbolically, denoted by $SN(\lambda)$

After the introduction of this distribution, a lot of research studies are carried out to introduce different skew normal distributions (see Chakraborty and Hazarika, 2011). Also, in the same spirit of the Azzalini skew normal distribution, Wahed and Ali (2001) and Nekoukhou and Alamatsaz (2012) proposed skew Logistic and skew Laplace distributions, respectively. In 2002, Balakrishnan introduced the generalization of SN($\lambda$) distribution with the pdf given by

$$f_Z(z;\lambda,n) = \varphi(z)[\Phi(\lambda z)]^n / C_n(\lambda); \quad -\infty < z < \infty, \lambda \in R \qquad (2)$$

where, $n$ is a positive integer and $C_n(\lambda) = E(\Phi^n(\lambda U))$, $U \sim N(0,1)$. For $\lambda = 1$, it reduces to SN($\lambda$) in equation (1). This distribution is known as Balakrishnan skew normal distribution. There are different extensions and generalizations of this distribution (For details see Sharafi and Behboodian (2008), Yadegari et al. (2008), Bahrami et al. (2009), Hasanalipour and Sharafi (2012) and among others).

Instead of using the cdf $\Phi(.)$ in equation (1), Huang and Chen (2007) introduced a new concept of skew function $G(.)$, a Lebesgue measurable function with $0 \le G(z) \le 1$ and $G(z) + G(-z) = 1$, $z \in R$, almost everywhere. The corresponding pdf is given by

$$f(z) = 2h(z)G(z) \; ; z \in R \qquad (3)$$

Employing this idea Elal-Olivero (2010) developed alpha skew normal distribution having pdf

$$f(z;\alpha) = \left(\frac{(1-\alpha z)^2 + 1}{2 + \alpha^2}\right)\varphi(z); \; z \in R, \alpha \in R \qquad (4)$$

The different generalization of the pdf in equation (4) were studied and analyzed (for details see Venegas et al. (2016), Louzada et al. (2017), Sharafi et al. (2017) etc.). A new form of skew logistic distribution and its generalization was given by Chakraborty et al. (2012) and Hazarika and Chakraborty (2015), respectively; and a multimodal form of skew Laplace distribution was studied by Chakraborty et al. (2014). A multimodal skewed extension of normal distribution was introduced by Chakraborty et al. (2015) by employing periodic skew function.

Using exactly the similar approach of Elal-Olivero (2010), Harandi and Alamatsaz (2013) and Hazarika and Chakraborty (2014) investigated a class of alpha skew Laplace distributions and alpha skew logistic distributions, respectively. In the present article, following the idea of Balakrishnan (2002) and taking $n = 2$, a new alpha skew logistic distribution is proposed and some of its basic properties are investigated.

## 2. New Alpha Skew Logistic Distribution

In this section we define a new alpha skew logistic distribution as follows

**Definition 1:** A continuous random variable $Z$ with the following pdf

$$f_Z(z;\alpha) = \frac{[(1-\alpha z)^2 + 1]^2}{C_2(\alpha)} \frac{e^{-z}}{(1 + e^{-z})^2}; \; z \in R \qquad (5)$$

is said to follow Balakrishnan alpha skew logistic distribution with skewness parameter $\alpha \in R$. $C_2(\alpha) > 0$ is the normalizing constant given in the Appendix D. In the rest of this article we shall refer the distribution in equation (5) as $BASLG_2(\alpha)$.

**Particular cases:**

- when $\alpha = 0$, we get the standard logistic distribution given by

$$f_Z(z) = \frac{e^{-z}}{(1+e^{-z})^2} \qquad (6)$$

and is denoted by $Z \to LG(0,1)$.

- when $\alpha \to \pm\infty$, we get a Bimodal Logistic ($BLG$) distribution (see Hazarika and Chakraborty, 2014) given by

$$f_Z(z) = \frac{15 z^4}{7 \pi^4} \frac{e^{-z}}{(1+e^{-z})^2} \qquad (7)$$

and is denoted by $Z \xrightarrow{d} BLG(4)$.

- If $Z \sim BASLG_2(\alpha)$ then $-Z \sim BASLG_2(-\alpha)$

The plot of the pdf of $BASLG_2(\alpha)$ is shown in figure 1 for studying the variation in its shape with respect to the parameter $\alpha$ clearly shows the role of $\alpha$ plays in skewness and kurtosis of the pdf. It's observed that the skewness is positive (negative) according as $\alpha >(<)0$ while for $\alpha = 0$ its symmetrical.

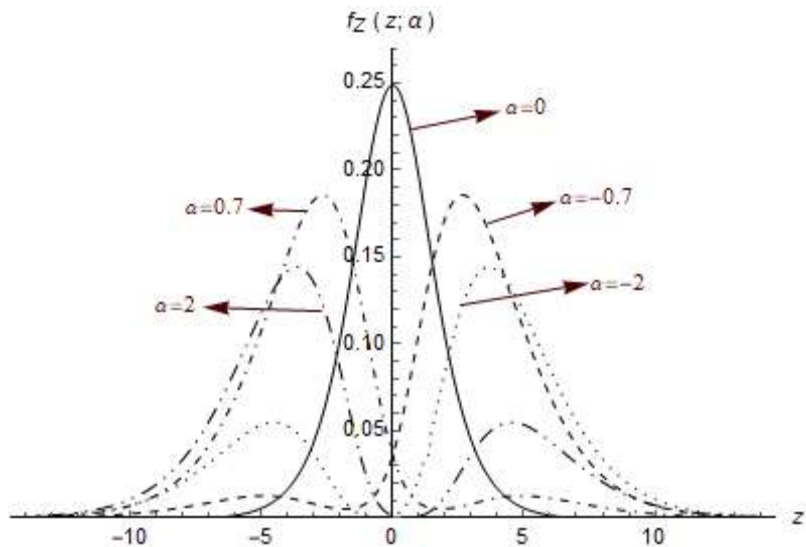

**Figure 1:** Plots of pdf of $BASLG_2(\alpha)$

**Theorem 1:** Let $Z \sim BASLG_2(\alpha)$ then its cdf is given by

$$F_Z(z;\alpha) = \frac{1}{(1+e^z)C_2(\alpha)} \begin{pmatrix} (2+z\alpha(-2+z\alpha))(e^z(2+z\alpha(-2+z\alpha))-4(1+e^z)\alpha(-1+z\alpha)Log[1+e^z]) - \\ 4(1+e^z)\alpha^2((4+3z\alpha(-2+z\alpha))Li_2[-e^z]+6\alpha((1-z\alpha)Li_3[-e^z]+\alpha Li_4[-e^z])) \end{pmatrix}$$

(8)

Where $Li_n[z] = \sum_{k=1}^{\infty} \frac{z^k}{k^n}$ is the poly-logarithm function (Prudnikov et al., 1986).

*Proof:* see Appendix A.

In particular when $\alpha \to \pm\infty$, $F_Z(z;\alpha)$ in equation (9) tends to the cdf of $BLG$ (4) distribution as

$$F_Z(z) = \frac{15}{7\pi^4}\left(\frac{e^z z^4}{1+e^z} - 4z^3 Log[1+e^z] - 12z^2 Li_2[-e^z] + 24z Li_3[-e^z] - 24 Li_4[-e^z]\right) \quad (9)$$

**Plots of cdf:**

The plots of the cdf of $BASLG_2(\alpha)$ in figure 2 exhibit the different variation in its shape with respect to the parameter $\alpha$.

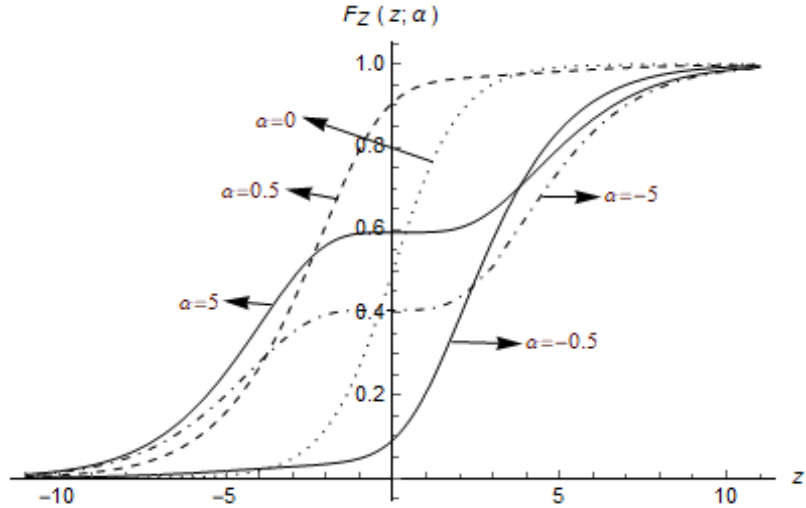

**Figure 2:** Plots of cdf of $BASLG_2(\alpha)$

From the above figure it can be easily seen that for the value of the parameter $\alpha \neq 0$, the tails get stretched, more over for larger value of $|\alpha|$ the bimodality is prominently expressed (see also figure 1).

**Theorem 2:** The $BASLG_2(\alpha)$ have at most two modes.

*Proof:* The first differentiation of the pdf of $BASLG_2(\alpha)$ with respect to z gives

$$Df_Z(z;\alpha) = \frac{e^z\left((1-\alpha z)^2 + 1\right)\left(-2 + 2\alpha(z+2) - z\alpha^2(z+4) + \alpha e^z\left(2 - 2(z-2) + \alpha^2 z(z-4)\right)\right)}{(1+e^z)^3 C_2(\alpha)} \quad (10)$$

Then the contour of the equation $Df_Z(z;\alpha) = 0$ is drawn in figure 3. From the figures it can be seen that there is at most three zeros of $Df_Z(z;\alpha)$ which implies that $BASLG_2(\alpha)$ distribution has at most two modes. From the figure it can also be be noted that for $-0.48 < \alpha < 0.48$, $BASLG_2(\alpha)$ remains unimodal (see on the right panel of figure 3).

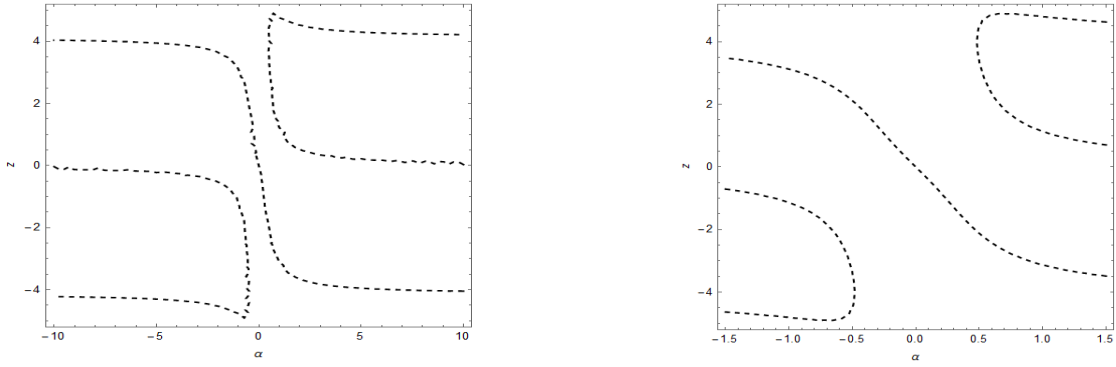

**Figure 3:** Contour plots of the equation $Df_Z(z;\alpha) = 0$

**Theorem 3:** Let $Z \sim BASLG_2(\alpha)$ then its moment generating function (mgf) is given by

$$M_Z(t) = \frac{1}{C_2(\alpha)} \begin{pmatrix} \pi\, Csc(\pi t)(-8\alpha + 12\pi^2\alpha^3 + t(4 - 8\pi^2\alpha^2 + \pi^4\alpha^4) \\ + 4\pi\alpha((-4\alpha + \pi^2\alpha^3 + t(2 - \pi^2\alpha^2))Cot(\pi t) + \pi\alpha(-6\alpha + t(4 - 5\pi^2\alpha^2) \\ + 6\pi(t-\alpha)\alpha\, Cot(\pi t))Csc(\pi t)^2 + 6\pi^3 t\alpha^3 Csc(\pi t)^4)) \end{pmatrix}$$

provided $-1 < t < 1$. (11)

*Proof:* see Appendix B.

In particular if $\alpha \to \pm\infty$, $M_Z(t)$ in equation (11) tends to the mgf of *BLG* (4) distribution as

$$M_Z(t) = \frac{15}{56} Csc(\pi t)^5 (115\pi t + \pi t(76 Cos(2\pi t) + Cos(4\pi t)) - 88 Sin(2\pi t) - 4 Sin(4\pi t)) \quad (12)$$

**Theorem 4:** Let $Z \sim BASLG_2(\alpha)$ then for $k \in N$, the $k^{th}$ moment is given by

$$E(Z^{2k}) = \frac{1}{C_2(\alpha)} 2\left[ 4\left(1 - \frac{1}{2^{k-1}}\right)\Gamma(k+1)\varsigma(k) + 8\alpha^2\left(1 - \frac{1}{2^{k+1}}\right)\Gamma(k+3)\varsigma(k+2) + \right.$$
$$\left. \alpha^4\left(1 - \frac{1}{2^{k+3}}\right)\Gamma(k+5)\varsigma(k+4) \right] \quad (13)$$

and,

$$E(Z^{2k-1}) = \frac{1}{C_2(\alpha)}\left[ 2\left(-8\alpha\left(1 - \frac{1}{2^k}\right)\Gamma(k+2)\varsigma(k+1) - 4\alpha^3\left(1 - \frac{1}{2^{k+2}}\right)\Gamma(k+5)\varsigma(k+3)\right) \right]$$

(14)

where, $\varsigma(s) = \sum_{j=1}^{\infty} j^{-s}$ is the Riemann zeta functions and $\Gamma(s) = \int_0^{\infty} t^{s-1} \exp(-t) dt$ (Prudnikov et al., 1986).

*Proof:* See Appendix C.

In particular it's easy to derive the following expressions for first four moments and variance:

$$E(Z) = \frac{4\pi^2\alpha(10 + 7\pi^2\alpha^2)}{60 + 40\pi^2\alpha^2 + 7\pi^4\alpha^4}$$

$$E(Z^2) = \frac{\pi^2(140 + 392\pi^2\alpha^2 + 155\pi^4\alpha^4)}{7(60 + 40\pi^2\alpha^2 + 7\pi^4\alpha^4)}$$

$$E(Z^3) = \frac{4\pi^4\alpha(98 + 155\pi^2\alpha^2)}{7(60 + 40\pi^2\alpha^2 + 7\pi^4\alpha^4)}$$

$$E(Z^4) = \frac{\pi^4(196 + 1240\pi^2\alpha^2 + 889\pi^4\alpha^4)}{7(60 + 40\pi^2\alpha^2 + 7\pi^4\alpha^4)}$$

$$Var(Z) = \frac{\pi^2(8400 + 17920\pi^2\alpha^2 + 10280\pi^4\alpha^4 + 3456\pi^6\alpha^6 + 1085\pi^8\alpha^8)}{7(60 + 40\pi^2\alpha^2 + 7\pi^4\alpha^4)}$$

The skewness and kurtosis of $BASLG_2(\alpha)$ can be easily obtained (given in the Appendix E).

**Remark 1:** By numerically optimizing mean, variance, skewness and the kurtosis with respect to $\alpha$, we obtain their bounds respectively as $-2.9077 \leq E(Z) \leq 2.9077$, $3.28987 \leq Var(Z) \leq 31.2202$, $0 \leq \beta_1 \leq 1.3945$ and $1.81315 \leq \beta_2 \leq 6.87571$.

The plots of mean, variance, skewness and kurtosis of $BASLG_2(\alpha)$ are shown, respectively, in figure 4, 5, 6 and 7 further validate the Remark 1.

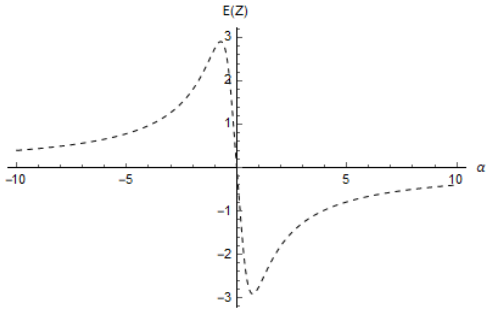

**Figure 4:** Plot of mean

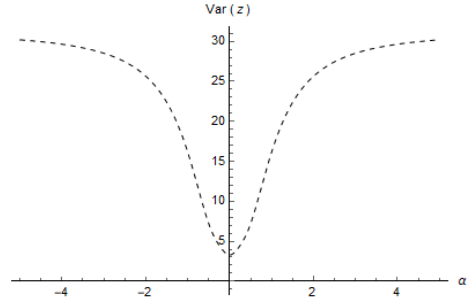

**Figure 5:** Plot of variance

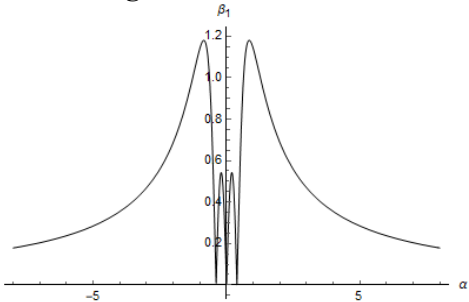

**Figure 6:** Plot of skewness

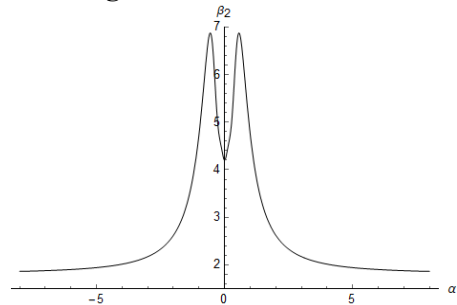

**Figure 7:** Plot of kurtosis

**Remark 2:** The pdf of equation (5) can be expressed as sum of two functions, as shown below

$$f_Z(z;\alpha) = \frac{[(1-\alpha z)^2 + 1]^2}{C_2(\alpha)} \frac{e^{-z}}{(1+e^{-z})^2}$$

$$= \frac{4 + 8\alpha^2 z^2 + \alpha^4 z^4}{C_2(\alpha)} \frac{e^{-z}}{(1+e^{-z})^2} - \frac{8\alpha z + 4\alpha^3 z^3}{C_2(\alpha)} \frac{e^{-z}}{(1+e^{-z})^2} \quad (15)$$

In equation (15), the first part is symmetric and the second part is asymmetric one and the symmetric part, which is defined below, is symbolically denoted by $SCBASLG_2(\alpha)$.

**Definition 2:** A continuous random variable $Z$ with the following pdf

$$f_1(z;\alpha) = \frac{4 + 8\alpha^2 z^2 + \alpha^4 z^4}{C_2(\alpha)} \frac{e^{-z}}{(1+e^{-z})^2}; \quad z \in R, \alpha \in R \tag{16}$$

is the symmetric-component random variable of the model $BASLG_2(\alpha)$. We write it as $Z \sim SCBASLG_2(\alpha)$.

**Properties of $SCBASLG_2(\alpha)$:**

(i) If $\alpha = 0$, then $Z \sim LG(0,1)$.

(ii) If $\alpha \to \pm\infty$, then $Z \xrightarrow{d} BLG$

**Theorem 5:** Let $Z \sim SCBASLG_2(\alpha)$ then its cdf is given by

$$F_1(z;\alpha) = \frac{1}{C_2(\alpha)(1+e^z)} \begin{pmatrix} e^z(4 + 8\alpha^2 z^2 + \alpha^4 z^4) - 4(1+e^z)\alpha^2 z(4+\alpha^2 z^2)Log[1+e^z] - \\ 4(1+e^z)\alpha^2((4+3\alpha^2 z^2)Li_2[-e^z] + 6\alpha^2(-zLi_3[-e^z] + Li_4[-e^z])) \end{pmatrix}$$

(17)

Proof: see Appendix F.

**Theorem 6:** Let $Z \sim SCBASLG_2(\alpha)$ then its mgf is given by

$$M_Z(t;\alpha) = \frac{\pi\, Csc[\pi t]^5}{8C_2(\alpha)(1+e^z)} \begin{pmatrix} 12t + 40\pi^2\alpha^2 t + 115\pi^4\alpha^4 t + 4tCos[2\pi t](-4 - 8\pi^2\alpha^2 + 19\pi^4\alpha^4) + \\ tCos[4\pi t](4 - 8\pi^2\alpha^2 + \pi^4\alpha^4) - 32\pi\alpha^2 Sin[2\pi t](1 + 2.75\pi^2\alpha^2) + \\ 4\pi\alpha^2 Sin[4\pi t](4 - \pi^2\alpha^2) \end{pmatrix}$$

provided, $-1 < t < 1$ (18)

Proof: see Appendix G.

**Remark 3:** The following algorithm can be used to simulate random numbers for $BASLG_2(\alpha)$ distribution:

Let $g(x;\alpha)$, the density function of $X \sim BASLG_2(\alpha)$ and $g_1(x;\alpha)$, the density function of $Y \sim SCBASLG_2(\alpha)$, with

$$S = \underset{x}{Sup}\, \frac{g(x;\alpha)}{g_1(x;\alpha)} = \frac{3 + 2\sqrt{2}}{3}$$

Then a random variable from $BASLG_2(\alpha)$ distribution can be generated in the following steps:

Let $X$ and $Y$ are two independent random variables then

(i) Generate $Y \sim SCBASLG_2(\alpha)$ and $U \sim Uniform(0,1)$, and they are independent.

(ii) If $U < \frac{1}{S}\frac{g(x;\alpha)}{g_1(x;\alpha)}$, and set $X = Y$; otherwise go back to step(i) and continue the process.

## 3. Some Extensions of $BASLG_2(\alpha)$

In this section we briefly mention some of the possible extensions of $BASLG_2(\alpha)$ distributions. These are being currently investigated and will be reported later.

## 3.1. Location and Scale Extension

If $Z \sim BASLG_2(\alpha)$ then for $\mu \in R$ and $\beta > 0$, set $Y = \beta Z + \mu$. Therefore, the location-scale extension of $BASLG_2(\alpha)$ has the following pdf given by

$$f_Y(y;\mu,\beta,\alpha) = \frac{\left[\left(1-\alpha\left(\frac{y-\mu}{\beta}\right)\right)^2 + 1\right]^2 e^{-\left(\frac{y-\mu}{\beta}\right)}}{\beta C_2(\alpha)\left(1+e^{-\frac{y-\mu}{\beta}}\right)^2} ; \quad y,\mu,\alpha \in R, \beta > 0 \qquad (19)$$

Symbolically, we write, $Y \sim BASLG_2(\alpha,\mu,\beta)$.

## 3.2. The Bivariate extension

In this section, using the idea of Louzada et al. (2017) a natural bivariate form of $BASLG_2(\alpha)$ distribution is is defined as follows:

**Definition 3:** A two-dimensional random variable $\mathbf{Z} = (Z_1, Z_2)$ with the following joint pdf

$$f_\mathbf{Z}(\mathbf{z};\alpha,\alpha_1,\alpha_2) = \frac{[(1-\alpha_1 z_1 - \alpha_2 z_2)^2 + 1]^2}{C_2(\alpha,\alpha_1,\alpha_2)} \Psi(z_1,z_2,\alpha) \qquad (20)$$

where, $C_2(\alpha,\alpha_1,\alpha_2)$ and $\Psi(z_1,z_2,\alpha)$ are given in Appendix D; $z_1, z_2, \alpha, \alpha_1, \alpha_2 \in R$. $\Psi(z_1,z_2,\alpha)$ is the pdf of bivariate logistic distribution (see Gumbel, 1961) is said to be Bivariate Balakrishnan alpha skew logistic distribution with the parameters $\alpha, \alpha_1$ and $\alpha_2$ and is denoted by $BBASLG_2(\alpha,\alpha_1,\alpha_2)$.

## 3.3. A Two-parameter $BASLG_2(\alpha)$

Like Bahrami et al. (2008), here we present the definition and some simple properties of a two-parameter Balakrishnan alpha skew normal distribution.

**Definition 4:** If a continuous random variable Z has the following pdf

$$f_Z(z;\alpha_1,\alpha_2) = \frac{[(1-\alpha_1 z)^2 + 1]^2 [(1-\alpha_2 z)^2 + 1]^2}{C_2(\alpha_1,\alpha_2)} \frac{e^{-z}}{(1+e^{-z})^2}; \quad z \in R \qquad (21)$$

where, $C_2(\alpha_1,\alpha_2)$ is given in the Appendix D, it is said to follow two-parameter Balakrishnan alpha skew logistic distribution with parameters $\alpha_1, \alpha_2 \in R$, and is denoted by $TPBASLG_2(\alpha_1,\alpha_2)$.

## 3.4. Balakrishnan Alpha-Beta Skew Logistic Distribution

In this section, following the study of Shafiei et al. (2016), we present the definition and some simple properties of Balakrishnan alpha-beta skew logistic distribution.

**Definition 5:** A continuous random variable Z with the following pdf

$$f_Z(z;\alpha) = \frac{[(1-\alpha z - \beta z^3)^2 + 1]^2}{C_2(\alpha,\beta)} \frac{e^{-z}}{(1+e^{-z})^2}; \quad z \in R \qquad (22)$$

where $C_2(\alpha, \beta)$ is given in the Appendix D then, it is said to follow Balakrishnan alpha-beta skew logistic distribution with parameters $\alpha, \beta \in R$, and is denoted by $BABSLG_2(\alpha, \beta)$.

### 3.5. The Log-Balakrishnan Alpha Skew Logistic Distribution

Here, we present the definition and some simple properties of log-Balakrishnan alpha skew logistic distribution similar to the study of Venegas et al. (2016).

Let $Z = e^Y$, then $Y = Log(Z)$, therefore, the density function of Z is defined as follows:

**Definition 6:** A continuous random variable Z with the following pdf

$$f_Z(z;\alpha) = \frac{[(1 - \alpha\, Log\,[z])^2 + 1]^2}{C_2(\alpha)} \frac{e^{-Log[z]}}{z(1 + e^{-Log[z]})^2}; \quad z > 0 \qquad (23)$$

where $C_2(\alpha)$ is given in the Appendix D then, it is said to follow log-Balakrishnan alpha skew logistic distribution with parameter $\alpha \in R$, and is denoted by $LBASLG_2(\alpha, \beta)$.

### 4. Data Modelling Applications

#### 4.1. Maximum Likelihood Estimation

Let $Y \sim BASLG_2(\alpha, \mu, \sigma)$ then for a random sample $y_1, y_2, ..., y_n$ of size $n$, the log likelihood function is given by

$$l(\theta) = 2\sum_{i=1}^{n} \log\left[\left\{1 - \alpha\left(\frac{y_i - \mu}{\beta}\right)\right\}^2 + 1\right] - n\log(60 + 40\pi^2\alpha^2 + 7\pi^4\alpha^4) - n\log[15] - n\log(\beta)$$

$$- \sum_{i=1}^{n}\frac{(y_i - \mu)}{\beta} - 2\sum_{i=1}^{n} Log\left[1 + Exp\left[\frac{-(y_i - \mu)}{\beta}\right]\right]$$

Where, $\theta = (\alpha, \mu, \sigma)$ \hfill (24)

The MLE's of the parameters of $BASLG_2(\alpha, \mu, \beta)$ distribution are estimated numerically by maximizing the likelihood function using GenSA package in R.

For the applicability of the proposed model, we analyze the model with three real life datasets obtained from different sources. We then compared our proposed distribution $BASLG_2(\alpha, \mu, \beta)$ with the normal distribution $N(\mu, \sigma^2)$, the logistic distribution $LG(\mu, \sigma)$, the Laplace distribution $La(\mu, \beta)$, the skew-normal distribution $SN(\lambda, \mu, \sigma)$ of Azzalini (1985), the skew-logistic distribution $SLG(\lambda, \mu, \beta)$ of Wahed and Ali (2001), the skew-Laplace distribution $SLa(\lambda, \mu, \beta)$ of Nekoukhou and Alamatsaz (2012), the alpha-skew-normal distribution $ASN(\alpha, \mu, \sigma)$ of Elal-Olivero (2010), the alpha-skew-Laplace distribution $ASLa(\alpha, \mu, \beta)$ of Harandi and Alamatsaz (2013), the alpha-skew-logistic distribution $ASLG(\alpha, \mu, \beta)$ of Hazarika and Chakraborty (2014), the alpha-beta-skew-normal distribution $ABSN(\alpha, \beta, \mu, \sigma)$ and the beta-skew-normal distribution $BSN(\beta, \mu, \sigma)$ of Shafiei et al. (2016).

The MLE's of the parameters of different distributions mentioned above are estimated using numerical optimization schedule. AIC and BIC are used for model selection or comparison.

**Example 1:** N latitude degrees dataset:

This dataset is related to N latitude degrees in 69 samples from world lakes, which appear in Column 5 of the Diversity data set in website: http://users.stat.umn.edu/sandy/courses/8061/datasets/lakes.lsp. Table 1 shows the MLE's, log-likelihood, AIC and BIC of the above mentioned distributions. The graphical representation of the results taking only the top three competitors for the proposed model is given in figure 8.

**Table 1:** MLE's, log-likelihood, AIC and BIC for N latitude degrees in 69 samples from world lakes.

| Parameters Distribution | $\mu$ | $\sigma$ | $\lambda$ | $\alpha$ | $\beta$ | $\log L$ | AIC | BIC |
|---|---|---|---|---|---|---|---|---|
| $N(\mu, \sigma^2)$ | 45.165 | 9.549 | -- | -- | -- | -253.599 | 511.198 | 515.666 |
| $LG(\mu, \beta)$ | 43.639 | -- | -- | -- | 4.493 | -246.65 | 497.290 | 501.758 |
| $SN(\lambda, \mu, \sigma)$ | 35.344 | 13.7 | 3.687 | -- | -- | -243.04 | 492.072 | 498.774 |
| $BSN(\beta, \mu, \sigma)$ | 54.47 | 5.52 | -- | -- | 0.74 | -242.53 | 491.060 | 497.760 |
| $SLG(\lambda, \mu, \beta)$ | 36.787 | -- | 2.828 | -- | 6.417 | -239.05 | 490.808 | 490.808 |
| $La(\mu, \beta)$ | 43 | -- | -- | -- | 5.895 | -239.25 | 482.496 | 486.964 |
| $ASLG(\alpha, \mu, \beta)$ | 49.087 | -- | -- | 0.861 | 3.449 | -237.35 | 480.702 | 487.404 |
| $SLa(\lambda, \mu, \beta)$ | 42.3 | -- | 0.255 | -- | 5.943 | -236.90 | 479.799 | 486.501 |
| $ASLa(\alpha, \mu, \beta)$ | 42.3 | -- | -- | -0.22 | 5.44 | -236.08 | 478.159 | 484.861 |
| $ASN(\alpha, \mu, \sigma)$ | 52.147 | 7.714 | -- | 2.042 | -- | -235.37 | 476.739 | 483.441 |
| $ABSN(\alpha, \beta, \mu, \sigma)$ | 53.28 | 9.772 | -- | 2.943 | -0.292 | -234.36 | 476.719 | 485.655 |
| $BASLG_2(\alpha, \mu, \beta)$ | 52.494 | -- | -- | 0.907 | 2.638 | **-230.75** | **467.501** | **474.203** |

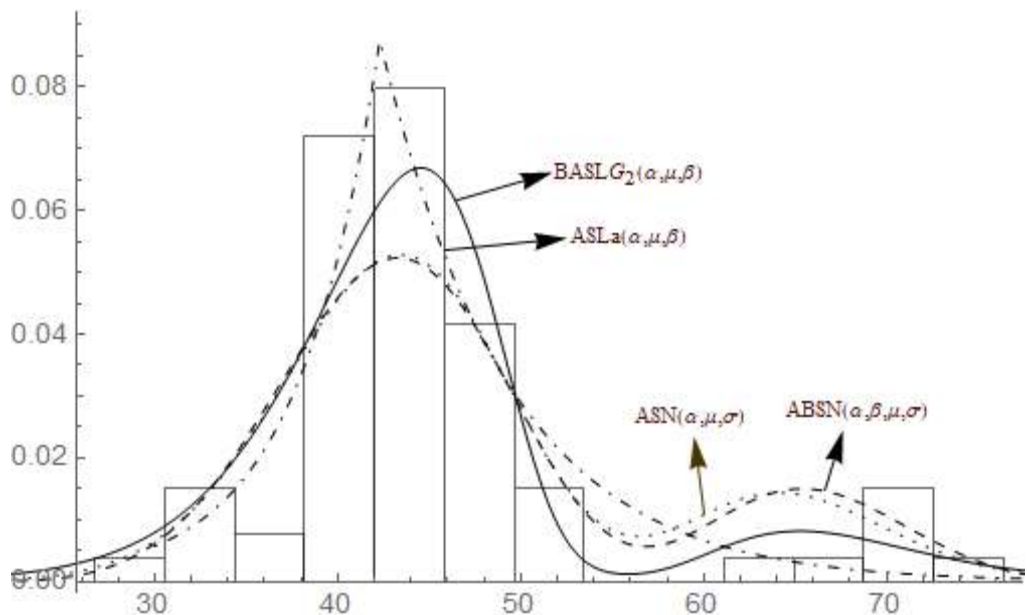

**Figure 8:** Plots of observed and expected densities of some distributions for N latitude degrees in 69 samples from world lakes.

**Example 2:** Exchange rate dataset:

To illustrate more we used the exchange rate data of the United Kingdom Pound to the United States Dollar from 1800 to 2003. The data obtained from the website http://www.globalfindata.com. The MLE's, log-likelihood, AIC and BIC of the above given distributions are shown in Table 2 and graphical representation are given in figure 9.

**Table 2:** MLE's, log-likelihood, AIC and BIC for the exchange rate data of the United Kingdom Pound to the United States Dollar from 1800 to 2003.

| Parameters<br>Distribution | $\mu$ | $\sigma$ | $\lambda$ | $\alpha$ | $\beta$ | $\log L$ | AIC | BIC |
|---|---|---|---|---|---|---|---|---|
| $N(\mu,\sigma^2)$ | 4.117 | 1.381 | -- | -- | -- | -355.265 | 714.529 | 721.165 |
| $SN(\lambda,\mu,\sigma)$ | 3.589 | 1.478 | 0.501 | -- | -- | -355.217 | 716.434 | 726.388 |
| $LG(\mu,\beta)$ | 4.251 | -- | -- | -- | 0.753 | -351.192 | 706.385 | 713.021 |
| $SLG(\lambda,\mu,\beta)$ | 5.36 | -- | -2.371 | -- | 1.018 | -341.391 | 688.782 | 698.736 |
| $La(\mu,\beta)$ | 4.754 | -- | -- | -- | 0.971 | -339.315 | 682.63 | 689.265 |
| $BSN(\beta,\mu,\sigma)$ | 4.526 | 0.974 | -- | -- | 0.181 | -334.488 | 674.975 | 684.93 |
| $ASN(\alpha,\mu,\sigma)$ | 3.656 | 0.883 | -- | -3.504 | -- | -317.946 | 641.892 | 651.847 |
| $SLa(\lambda,\mu,\beta)$ | 4.855 | -- | 1.506 | -- | 1 | -311.318 | 628.636 | 638.59 |
| $ABSN(\alpha,\beta,\mu,\sigma)$ | 3.581 | 1.122 | -- | -5.542 | 0.715 | -301.788 | 611.576 | 624.848 |
| $ASLG(\alpha,\mu,\beta)$ | 3.764 | -- | -- | -2.025 | 0.403 | -301.963 | 609.927 | 619.881 |
| $ASLa(\alpha,\mu,\beta)$ | 4.861 | -- | -- | 0.539 | 0.677 | -301.443 | 608.885 | 618.84 |
| $BASLG_2(\alpha,\mu,\beta)$ | 3.671 | -- | -- | -1.693 | 0.265 | -300.583 | **607.166** | **617.121** |

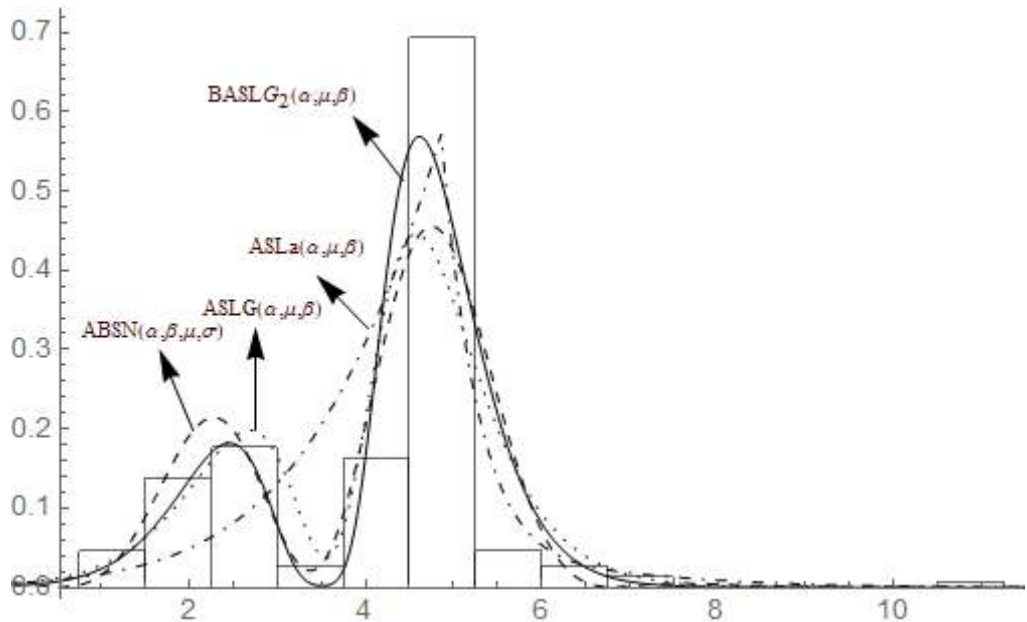

**Figure 9:** Plots of observed and expected densities of some distributions for the exchange rate data of the United Kingdom Pound to the United States Dollar from 1800 to 2003.

**Example 3:** Velocities of 82 distant galaxies dataset:

The third dataset consists of the velocities of 82 distant galaxies, diverging from our own galaxy. The data set is available at http://www.stats.bris.ac.uk/~peter/mixdata. The results for the above mentioned distributions are shown in Table 3 and their graphical representations are given in figure 10.

**Table 3:** MLE's, log-likelihood, AIC and BIC for the velocities of 82 distant galaxies, diverging from our own galaxy.

| Parameters Distribution | $\mu$ | $\sigma$ | $\lambda$ | $\alpha$ | $\beta$ | $\log L$ | AIC | BIC |
|---|---|---|---|---|---|---|---|---|
| $N(\mu, \sigma^2)$ | 20.832 | 4.54 | -- | -- | -- | -240.42 | 484.833 | 489.646 |
| $SN(\lambda, \mu, \sigma)$ | 24.61 | 5.907 | -1.395 | -- | -- | -239.21 | 484.42 | 491.64 |
| $SLG(\lambda, \mu, \beta)$ | 21.532 | -- | -0.154 | -- | 2.219 | -233.31 | 472.628 | 479.849 |
| $LG(\mu, \beta)$ | 21.075 | -- | -- | -- | 2.204 | -233.65 | 471.299 | 476.113 |
| $BSN(\beta, \mu, \sigma)$ | 20.596 | 3.26 | -- | -- | -0.158 | -232.22 | 470.44 | 477.66 |
| $ASN(\alpha, \mu, \sigma)$ | 17.417 | 3.869 | -- | -1.656 | -- | -230.09 | 466.175 | 473.395 |
| $SLa(\lambda, \mu, \beta)$ | 20.846 | -- | 1.002 | -- | 2.997 | -228.83 | 463.658 | 470.878 |
| $La(\mu, \beta)$ | 20.838 | -- | -- | -- | 2.997 | -228.83 | 461.66 | 466.474 |
| $ASLG(\alpha, \mu, \beta)$ | 18.482 | -- | -- | -0.833 | 1.646 | -224.88 | 455.754 | 462.974 |
| $ABSN(\alpha, \beta, \mu, \sigma)$ | 19.448 | 3.462 | -- | -1.392 | 0.323 | -220.06 | 448.109 | 457.736 |
| $ASLa(\alpha, \mu, \beta)$ | 19.473 | -- | -- | -0.842 | 1.805 | -220.79 | 447.586 | 454.806 |
| $BASLG_2(\alpha, \mu, \beta)$ | 17.117 | -- | -- | -0.799 | 1.263 | -219.86 | **445.716** | **452.936** |

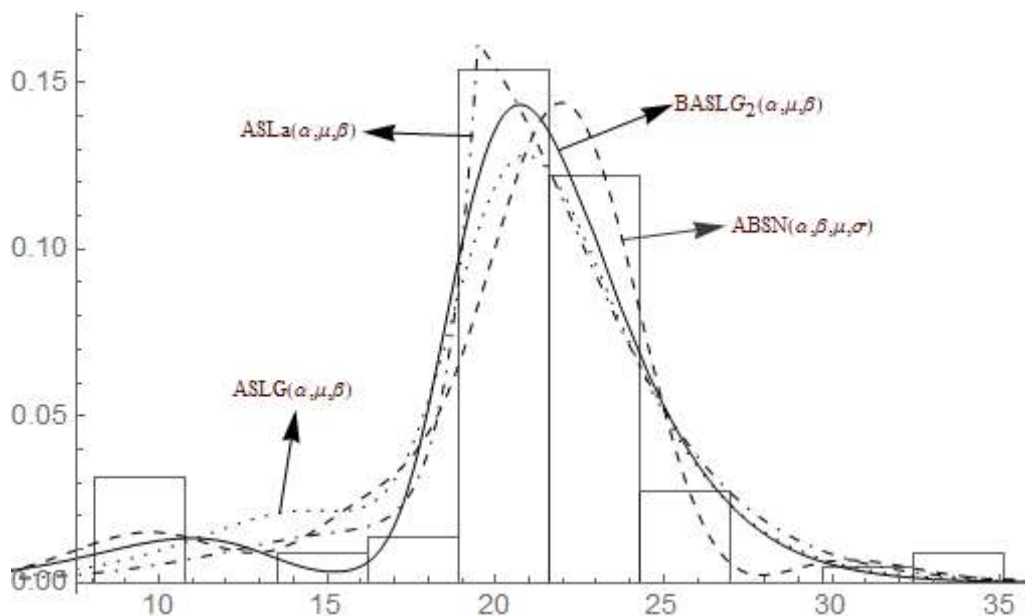

**Figure 10:** Plots of observed and expected densities of some distributions for the velocities of 82 distant galaxies, diverging from our own galaxy.

## 4.2. Likelihood Ratio Test

The likelihood ratio (LR) test is used to discriminate between $LG(\mu, \beta)$ and $BASLG_2(\alpha, \mu, \beta)$ as they are nested models. Here the Null hypothesis $H_0 : \alpha = 0$, that is the sample is drawn from $LG(\mu, \beta)$; against the alternative $H_1 : \alpha \neq 0$, that is the sample is drawn from $BASLG_2(\alpha, \mu, \beta)$.

The values of the LR test statistic for the above three examples are given in Table 4.

**Table 4:** the values of the LR test statistic for different datasets

| Datasets | Example 1 | Example 2 | Example 3 |
|---|---|---|---|
| LR | 31.7894 | 101.2184 | 27.5838 |

Since the calculated value of LR test statistic for all the datasets are greater than 6.635 the 99% tabulated value of chi-square at 1 degree of freedom, therefore there is no evidence in support of the null hypothesis for all the datasets. Thus, we conclude that for all the examples the sampled data comes from $BASLG_2(\alpha, \mu, \beta)$, not from $LG(\mu, \beta)$.

## 5. Conclusions

In this study a new bimodal alpha-skew-logistic distribution with one parameter is constructed and its properties are studied. The proposed distribution has nice compact form for cdf and moments and has provided best fitting to the data sets considered here in terms of AIC and BIC. It is therefore expected that the new distribution will be a useful one for data modelling.

# Appendix

**A: Proof of Theorem 1**

$$F_Z(z;\alpha) = \frac{1}{C_2(\alpha)} \int_{-\infty}^{z} [(1-\alpha z)^2 + 1]^2 \frac{e^{-z}}{(1+e^{-z})^2} dz$$

$$= \frac{1}{C_2(\alpha)} \left[ \int_{-\infty}^{z} \frac{4e^{-z}}{(1+e^{-z})^2} dz - 8\alpha \int_{-\infty}^{z} \frac{z e^{-z}}{(1+e^{-z})^2} dz + 8\alpha^2 \int_{-\infty}^{z} \frac{z^2 e^{-z}}{(1+e^{-z})^2} dz - 4\alpha^3 \int_{-\infty}^{z} \frac{z^3 e^{-z}}{(1+e^{-z})^2} dz \right.$$

$$\left. + \alpha^4 \int_{-\infty}^{z} \frac{z^4 e^{-z}}{(1+e^{-z})^2} dz \right]$$

$$= \frac{1}{C_2(\alpha)} \left[ \frac{4e^z}{1+e^z} - \frac{8\alpha z e^z}{1+e^z} + 8\alpha^2 z \left( \left( \frac{z e^z}{1+e^z} - 2 \, Log[1+e^z] \right) - 2 \, Li_2[-e^z] \right) - \right.$$

$$4\alpha^3 \left( z^2 \left( z - \frac{z}{1+e^z} - 3 \, Log[1+e^z] \right) - 6z \, Li_2[-e^z] + 6 \, Li_3[-e^z] \right) +$$

$$\left. \alpha^4 \left( \frac{e^z z^4}{1+e^z} - 4z^3 Log[1+e^z] - 12 z^2 Li_2[-e^z] + 24 z Li_3[-e^z] - 24 Li_4[-e^z] \right) \right]$$

On Simplification we get the result in equation (8).

**B: Proof of Theorem 3**

$$M_Z(t) = \frac{1}{C_2(\alpha)} \int_{-\infty}^{\infty} e^{tz} \frac{[(1-\alpha z)^2 + 1]^2 e^{-z}}{(1+e^{-z})^2} dz$$

$$= \frac{1}{C_2(\alpha)} \left[ \int_{-\infty}^{\infty} \frac{4 e^{z(t-1)}}{(1+e^{-z})^2} dz - 8\alpha \int_{-\infty}^{\infty} \frac{z e^{z(t-1)}}{(1+e^{-z})^2} dz + 8\alpha^2 \int_{-\infty}^{\infty} \frac{z^2 e^{z(t-1)}}{(1+e^{-z})^2} dz - 4\alpha^3 \int_{-\infty}^{\infty} \frac{z^3 e^{z(t-1)}}{(1+e^{-z})^2} dz \right.$$

$$\left. + \alpha^4 \int_{-\infty}^{\infty} \frac{z^4 e^{z(t-1)}}{(1+e^{-z})^2} dz \right]$$

$$= \frac{1}{C_2(\alpha)} \left[ 4(\pi t \, Csc(\pi t)) - 8\alpha(\pi(1 - \pi t \, Cot(\pi t))Csc(\pi t)) + 8\alpha^2 \left( \frac{1}{2}\pi^2 \, Csc(\pi t)^3 (\pi t(3 + Cos(2\pi t)) - \right.\right.$$

$$2 Sin(2\pi t))) - 4\alpha^3 \left( -\frac{1}{4}\pi^3 \, Csc(\pi t)^4 (\pi t(23 \, Cos(\pi t) + Cos(3\pi t)) - 3(5 \, Sin(\pi t) + Sin(3\pi t))) \right) +$$

$$\left. \alpha^4 \left( \frac{1}{8}\pi^4 \, Csc(\pi t)^5 (115\pi t + \pi t(76 \, Cos(2\pi t) + Cos(4\pi t)) - 88 \, Sin(2\pi t) + Sin(4\pi t)) \right) \right]$$

On simplification we get the desired result in equation (11).

**C: Proof of Theorem 4**

Since, the $k^{th}$ moment of the standard logistic distribution (see Balakrishnan (1992)) is given by

$$E(X^k) = \begin{cases} 2\left(1 - \dfrac{1}{2^{k-1}}\right)\Gamma(k+1)\varsigma(k) \text{ ; when } k \text{ is even} \\ 0 \text{ ; when } k \text{ is odd} \end{cases} \quad (A1)$$

Now, $E(Z^k) = \dfrac{1}{C_2(\alpha)} \int_{-\infty}^{\infty} Z^k \dfrac{[(1-\alpha z)^2 + 1]^2 e^{-z}}{(1+e^{-z})^2} dz$

$$= \dfrac{1}{C_2(\alpha)} \left[ \int_{-\infty}^{\infty} \dfrac{4z^k e^{-z}}{(1+e^{-z})^2} dz - 8\alpha \int_{-\infty}^{\infty} \dfrac{z^{k+1} e^{-z}}{(1+e^{-z})^2} dz + 8\alpha^2 \int_{-\infty}^{\infty} \dfrac{z^{k+2} e^{-z}}{(1+e^{-z})^2} dz - 4\alpha^3 \int_{-\infty}^{\infty} \dfrac{z^{k+3} e^{-z}}{(1+e^{-z})^2} dz \right.$$
$$\left. + \alpha^4 \int_{-\infty}^{\infty} \dfrac{z^{k+4} e^{-z}}{(1+e^{-z})^2} dz \right]$$

**Case1:** Let $k$ is even, then we get

$$E(Z^k) = \dfrac{1}{C_2(\alpha)} \left( \int_{-\infty}^{\infty} \dfrac{4z^k e^{-z}}{(1+e^{-z})^2} dz + 8\alpha^2 \int_{-\infty}^{\infty} \dfrac{z^{k+2} e^{-z}}{(1+e^{-z})^2} dz + \alpha^4 \int_{-\infty}^{\infty} \dfrac{z^{k+4} e^{-z}}{(1+e^{-z})^2} dz \right)$$

$$= \dfrac{1}{C_2(\alpha)} \left[ 4\left(1 - \dfrac{1}{2^{k-1}}\right)\Gamma(k+1)\varsigma(k) + 16\alpha^2 \left(1 - \dfrac{1}{2^{k+1}}\right)\Gamma(k+3)\varsigma(k+2) + 2\alpha^4 \left(1 - \dfrac{1}{2^{k+3}}\right)\Gamma(k+5)\varsigma(k+4) \right] \quad \text{[Using eqn. (A1)]}$$

On simplification we get the result in equation (13).

**Case2:** Let $k$ is odd, then we get

$$E(Z^k) = \dfrac{1}{C_2(\alpha)} \left( -8\alpha \int_{-\infty}^{\infty} \dfrac{z^{k+1} e^{-z}}{(1+e^{-z})^2} dz - 4\alpha^3 \int_{-\infty}^{\infty} \dfrac{z^{k+3} e^{-z}}{(1+e^{-z})^2} dz \right)$$

$$= \dfrac{1}{C_2(\alpha)} \left[ -16\alpha \left(1 - \dfrac{1}{2^k}\right)\Gamma(k+2)\varsigma(k+1) - 8\alpha^3 \left(1 - \dfrac{1}{2^{k+2}}\right)\Gamma(k+5)\varsigma(k+3) \right] \quad \text{[Using eqn. (A1)]}$$

On simplification we get the result (14).

**D: Normalizing Constants**

$$C_2(\alpha) = \left( 4 + \dfrac{8\pi^2 \alpha^2}{3} + \dfrac{7\pi^4 \alpha^4}{15} \right)$$

$$C_2(\alpha, \alpha_1, \alpha_2) = \dfrac{1}{15} \begin{pmatrix} 60 + 7\pi^4 \alpha_1^4 + 60\pi^2 \alpha \alpha_1^3 \alpha_2 + 40\pi^2 \alpha_2^2 + 7\pi^4 \alpha_2^4 + 10\pi^2 \alpha_1^2 (4 + \pi^2 \alpha_2^2) + \\ 60\alpha\alpha_1\alpha_2(4 + \pi^2 \alpha_2^2) \end{pmatrix}$$

$$\Psi(z_1, z_2, \alpha) = \dfrac{e^{-z_1 - z_2}}{(1+e^{-z_1})^2 (1+e^{-z_2})^2} \left[ 1 + \alpha\left(\dfrac{1-e^{-z_1}}{1+e^{-z_1}}\right)\left(\dfrac{1-e^{-z_2}}{1+e^{-z_2}}\right) \right]$$

$$C_2(\alpha_1, \alpha_2) = \dfrac{1}{105} \begin{pmatrix} 1680 + \pi^2(224\alpha_1\alpha_2(10 + 7\pi^2\alpha_2^2) + 28\alpha_2^2(40 + 7\pi^2\alpha_2^2) + 16\pi^2\alpha_1^3\alpha_2(98 + 155\pi^2\alpha_2^2) + \\ 8\alpha_1^2(140 + 392\pi^2\alpha_2^2 + 155\pi^4\alpha_2^4) + \pi^2\alpha_1^4(196 + 1240\pi^2\alpha_2^2 + 889\pi^4\alpha_2^4)) \end{pmatrix}$$

$$C_2(\alpha, \beta) = \dfrac{1}{15015} \begin{pmatrix} 60060 + 40040\pi^2\alpha^2 + 7007\pi^4\alpha^4 + 112112\pi^4\alpha\beta + 88660\pi^6\alpha^3\beta + 177320\pi^6\beta^2 + \\ 762762\pi^8\alpha^2\beta^2 + 465010\pi^{10}\alpha\beta^3 + 15559247\pi^{12}\beta^4 \end{pmatrix}$$

**E: Skewness and Kurtosis**

$$\beta_1 = \frac{11200\,\pi^2\alpha^2(10080 + 20240\,\pi^2\alpha^2 + 6320\,\pi^4\alpha^4 - 5112\,\pi^6\alpha^6 - 5194\,\pi^8\alpha^8 - 1519\,\pi^{10}\alpha^{10})^2}{(8400 + 17920\,\pi^2\alpha^2 + 10280\,\pi^4\alpha^4 + 3456\,\pi^6\alpha^6 + 1085\,\pi^8\alpha^8)^3}$$

$$\beta_2 = \frac{\begin{pmatrix} 7(42336000 + 207360000\,\pi^2\alpha^2 + 321427200\,\pi^4\alpha^4 + 306867200\,\pi^6\alpha^6 \\ + 231988640\,\pi^8\alpha^8 + 126595840\,\pi^{10}\alpha^{10} + 42918512\,\pi^{12}\alpha^{12} + 7353920\,\pi^{14}\alpha^{14} \\ + 304927\,\pi^{16}\alpha^{16}) \end{pmatrix}}{(8400 + 17920\,\pi^2\alpha^2 + 10280\,\pi^4\alpha^4 + 3456\,\pi^6\alpha^6 + 1085\,\pi^8\alpha^8)^2}$$

**F: Proof of Theorem 5**

$$F_1(z;\alpha) = \frac{1}{C_2(\alpha)} \int_{-\infty}^{z} \frac{(4 + 8\alpha^2 z^2 + \alpha^4 z^4)e^{-z}}{(1 + e^{-z})^2}\,dz$$

$$= \frac{1}{C_2(\alpha)} \left[ \int_{-\infty}^{z} \frac{4 e^{-z}}{(1 + e^{-z})^2}\,dz + 8\alpha^2 \int_{-\infty}^{z} \frac{z^2 e^{-z}}{(1 + e^{-z})^2}\,dz + \alpha^4 \int_{-\infty}^{z} \frac{z^4 e^{-z}}{(1 + e^{-z})^2}\,dz \right]$$

$$= \frac{1}{C_2(\alpha)} \left[ \frac{4 e^z}{1 + e^z} + 8\alpha^2 z\left( \left(\frac{z e^z}{1 + e^z} - 2\,Log[1 + e^z]\right) - 2\,Li_2[-e^z] \right) + \right.$$
$$\left. \alpha^4 \left( \frac{e^z z^4}{1 + e^z} - 4z^3\,Log[1 + e^z] - 12 z^2 Li_2[-e^z] + 24 z Li_3[-e^z] - 24 Li_4[-e^z] \right) \right]$$

On simplification we get the result in equation (17).

**G: Proof of Theorem 6**

$$M_Z(t;\alpha) = \frac{1}{C_2(\alpha)} \int_{-\infty}^{\infty} e^{tz}\frac{(4 + 8\alpha^2 z^2 + \alpha^4 z^4)e^{-z}}{(1 + e^{-z})^2}\,dz$$

$$= \frac{1}{C_2(\alpha)} \left[ \int_{-\infty}^{\infty} \frac{4 e^{z(t-1)}}{(1 + e^{-z})^2}\,dz + 8\alpha^2 \int_{-\infty}^{\infty} \frac{z^2 e^{z(t-1)}}{(1 + e^{-z})^2}\,dz + \alpha^4 \int_{-\infty}^{\infty} \frac{z^4 e^{z(t-1)}}{(1 + e^{-z})^2}\,dz \right]$$

$$= \frac{1}{C_2(\alpha)} \left[ 4(\pi t\,Csc[\pi t]) + 8\alpha^2 \left( \frac{1}{2}\pi^2\,Csc[\pi t]^3(\pi t(3 + Cos[2\pi t]) - 2\,Sin[2\pi t]) \right) + \right.$$
$$\left. \alpha^4 \left( \frac{1}{8}\pi^4\,Csc[\pi t]^5(115\,\pi t + \pi t(76\,Cos[2\pi t] + Cos[4\pi t]) - 88\,Sin[2\pi t] + Sin[4\pi t]) \right) \right]$$

On simplification we get the result in equation (18).